\documentclass{amsart}
\usepackage{verbatim,graphicx}

\newtheorem{theorem}{Theorem}[section]

\theoremstyle{definition}

\theoremstyle{remark}
\newtheorem{remark}[theorem]{Remark}

\numberwithin{equation}{section}

\usepackage[normalem]{ulem}

\usepackage[textsize=small]{todonotes}

\begin{document}

\title{Applications of topological graph theory to $2$-manifold learning}


\author{T. Berry}
\address{Department of Mathematical Sciences, George Mason University, 4400 University Drive, MS 3F2, Fairfax, VA 22030}
\email{tberry@gmu.edu}

\author{S. Schluchter}
\address{Department of Mathematical Sciences, George Mason University, 4400 University Drive, MS 3F2, Fairfax, VA 22030}
\email{steven.schluchter@gmail.com}



\begin{abstract}
We show how, given a sufficiently large point cloud sampled from an embedded 2-manifold in $\mathbb{R}^n$, we may obtain a global representation as a cell complex with vertices given by a representative subset of the point cloud.  The vertex spacing is based on obtaining an approximation of the tangent plane which insures that the vertex accurately summarizes the local data.  Using results from topological graph theory, we couple our cell complex representation with the known Classification of Surfaces in order to classify the manifold.  The algorithm developed gives a meaningful description of the embedding as a piecewise linear structure, which is obtained from combinatorial data by projecting points in the point cloud into estimates of tangent planes.\end{abstract}

\maketitle

\section{Introduction}

	Given a point cloud $X$ noisily sampled from an embedded manifold $f(M)$ in $\mathbb{R}^n$, the task is to discern information about $M$ or its embedding $f$. One existing method is based on the convergence of certain graph Laplacians to the intrinsic Laplacian operator on the manifold \cite{belkin2003laplacian,diffusion,SingerEstimate}, and while this method does capture global information about the manifold, it is extremely challenging to extract topological information such as simplicial homology.  The alternative is persistent homology \cite{ghrist2008,edels2010,carlsson2009topology}, which is explicitly targeted at simplicial homology, yet requires converting the dataset into an (abstract) simplicial complex which is an very large data structure.  Moreover, persistent homology does not return a particular homology, but rather many different homologies based on different scales.  While this may be a strength of the persistent homology method for some applications, in other contexts we may be interested in obtaining a representation of $f(M)$ from $X$.  Moreover, a significant challenge of persistent homology is obtaining an efficient subset of a large dataset for which the algorithm can be practically applied, a problem which will be addressed below.
The purpose of this article is to introduce tools from topological graph theory (the study of graphs in topological spaces) and recover a combinatorial approximation of $M$ in terms of a rotation system, which is a finite list of finite lists that encodes a $2$-complex whose $0$-skeleton is drawn from $X$.  Our Matlab implementation of this algorithm is available at \cite{Github}.

	Some distinct features of our methods are that the faces of a 2-complex $K$ encoded by a rotation system need not be triangles since a rotation system encodes a more general $2$-complex, and that we use the invariance of Euler characteristic under homotopy to accelerate the classification of $f(M)$ by contracting a spanning tree of the $1$-skeleton of $K$.
	
	Informally, our algorithm can be described as follows.  Given an aforementioned point cloud $X$, we pick any point $x\in X$, and we apply techniques from \cite{Berry} to find a good approximation for a plane $P$ tangent to $f(M)$ at $x$; we note the choice of a normal vector $v$ to $P$.  Points in $X$ deemed to be too close to both $P$ and $x$ to be of any significance are discarded.  We then orthogonally project the remaining points near $x$ into $P$ along $v$.  Looking down at $x$ in $P$ along $v$, we can now find a counterclockwise cyclic ordering, of the projected points, as they relate to $x$.  (One may use a convention of clockwise ordering, too.)  We then draw straight edges joining $x$ to the vertices that we kept, and we encode these edges and the order that they appear in terms of a list.  (This cyclic ordering of vertices neighboring a given a vertex in an embedded graph is called a rotation on $x$.)  We then look at the neighbors of $x$ from the vantagepoint induced by extending the local orientation we got from $v$, and we repeat this process until we have considered all of the points of $X$.  If at any point we traverse a cycle $C$ in the embedded graph, which reverses the orientation that we previously had by extending the normal vector $v$, then we declare the last edge $e$ of $C$ to be ``orientation reversing" (a property that we will explain in Section \ref{section:graphBackground}), and we note this property in the rotations on the vertices to which $e$ is incident.  
	
	The result is that for each vertex chosen from the data set we have a list of neighboring vertices ordered according to their angle when projected into the tangent space.  This finite list of finite lists is called a rotation system, and it encodes a piecewise linear approximation of $f(M)$ in the form of a $2$-complex $K$.  Using the Euler characteristic and orientability of $K$, and the known Classification of Surfaces \cite[Theorem 77.5]{M}, we may classify $K$ and declare $K$ to be a reasonable approximation of $f(M)$.
	
In Section \ref{section:graphBackground}, we state all requisite background from graph theory and topological graph theory.  In Section \ref{section:manifoldBackground}, we state and describe the theory we need from the study of manifold learning.  In Section \ref{section:theAlgorithm}, we describe our algorithm in detail, including all of the choices that were made when implementing our algorithm.  In Section \ref{section:data}, we describe our experimental results and give some closing remarks, including some suggestions for further investigation.

\section{Background information in graph theory and topological graph theory}\label{section:graphBackground}

For additional background material, the reader is encouraged to consult \cite[Chapters 1, 2 and 4]{GT}. 

In this article, a graph $G=(V,E)$ is a finite and connected multigraph, allowing for loops and parallel edges.  An edge is a link if it is not a loop.  A surface is a compact and connected $2$-manifold, without boundary.  A cellular embedding of a graph in a surface is an embedding having the property that the complement of the graph is homeomorphic to a disjoint union of discs.  We will let $G$, $S$, and $G\rightarrow S$ denote a graph, a surface, and a cellular embedding of $G$ in $S$, respectively.

\begin{figure}[h]
\begin{center}
\includegraphics[scale=0.5]{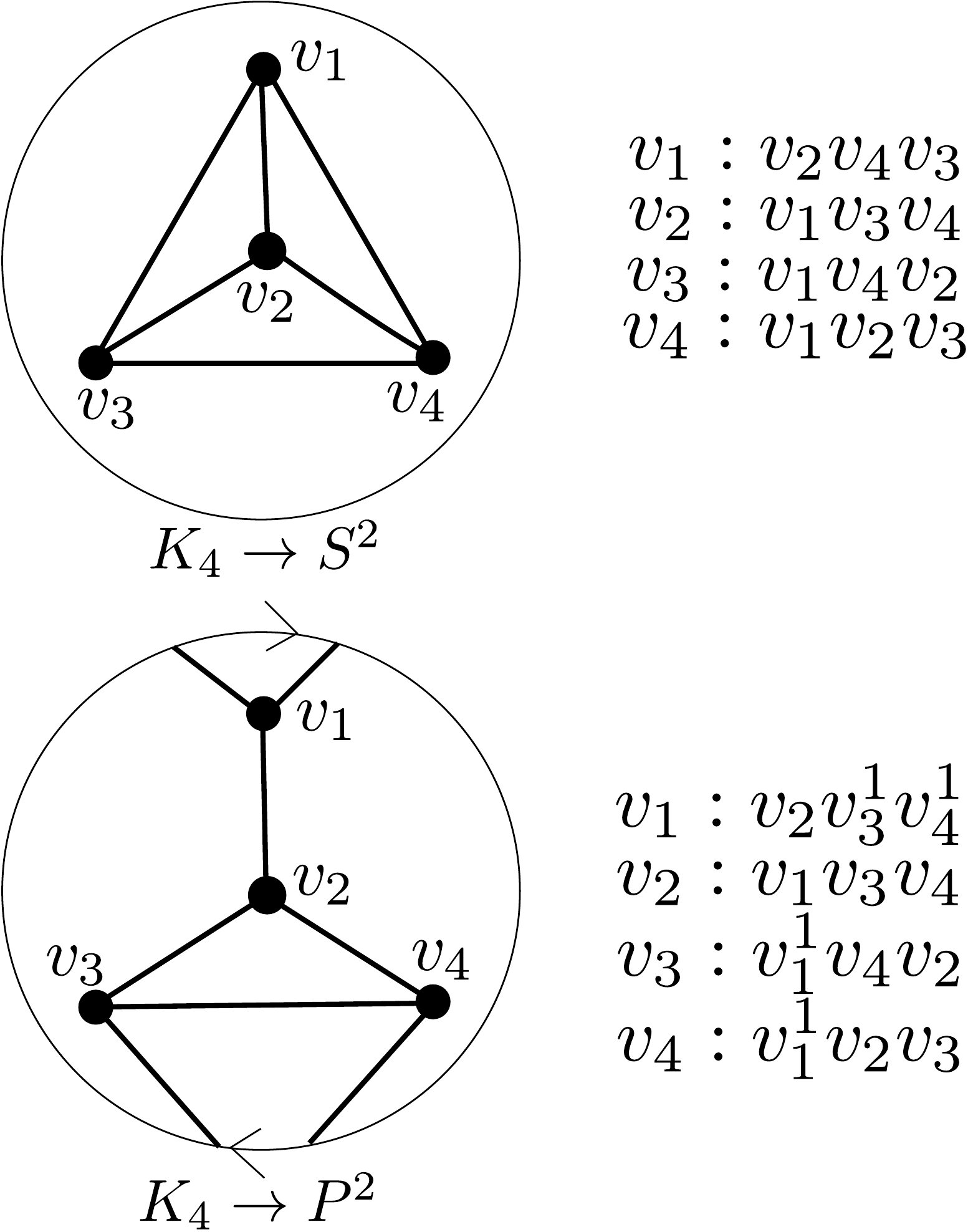}
\caption{\label{fig:rotationSystems} Rotation systems of orientable and nonorientable embeddings of $K_4$ in the sphere and the projective plane, respectively.  The superscrpts in the nonorientable embedding denote type-1 edges.}
\end{center}
\end{figure}

Given $G\rightarrow S$, if one thickens $G$ and then deletes the faces of $G\rightarrow S$, one obtains a band decomposition of $G\rightarrow S$: each vertex becomes a disc (a $0$-cell) and each edge becomes homeomorphic to a strip of ribbon (a $1$-cell), with the added condition that some of the strips may be twisted, reflecting the orientability of $S$.  We will make use of the notion of a surface with boundary in Section \ref{section:data}, and the only surface with boundary we will discuss will be the M\"obius band.  A twisted edge is called a type-1 edge, whereas a nontwisted edge is called a type-0 edge.  The union of the $0$-bands and $1$-bands corresponding to $G\rightarrow S$ forms what is called a ribbon graph of $G\rightarrow S$.  We may now give a combinatorial encoding of $G\rightarrow S$ exists in the form of what is called a rotation system, which can be thought of as an encoding of the ribbon graph of $G\rightarrow S$, which lists a cyclic ordering of the edges, with their edge types, incident to each vertex of $G$ (given a (counter)clockwise choice of ordering convention).  Examples of rotation systems corresponding to cellular embeddings of $K_4$ appear in Figure \ref{fig:rotationSystems}, where since the graphs are simple, and since each edge joins a unique pair of vertices, we may list the cyclic orderings of the adjacent vertices.

While each cellular graph embedding in a surface can be encoded in the form of a rotation system, it is not the case that the encodings are necessarily unique.  A local sign switch can be applied to a vertex of a rotation system, which is a reversal of the cyclic ordering and changing of the types of the links incident to that vertex.  (Local sign switches are used extensively in topological graph theory.  See \cite{Archdeacon} for more.)  The reader may benefit from thinking of a local sign switch of a vertex of $G\rightarrow S$ as turning over the corresponding $0$ band.  The reader should note that while a rotation system may feature several type-$1$ edges, the corresponding surface may be orientable: If each cycle of $G\rightarrow S$ contains an even number of type-$1$ edges, then all of the edges can be made type-$0$ edges by using local sign switches and $S$ made explicitly orientable.  See \cite[Example 3.2.2]{GT} for an example where the type-$1$ edges can all be ``untwisted", resulting in an orientable manifold.

Given $G\rightarrow S$ encoded in the form of a rotation system, it is possible to determine the classification of $S$.  We may use the known Classification of Surfaces \cite[Theorem 77.5]{M}, which states that $S$ is homeomorphic to: a sphere, the connected sum of a definitive number of torii, or the connected sum of a definitive number of projective planes.  All that we will need to determine the class of $S$ is the Euler characteristic and orientability of $S$.  The orientability of $S$ can be found by using local sign switches to determine if there is a spanning tree of $G$ which consists only of type-$0$ edges.  Since one may easily obtain the number of vertices and edges of $G$ from a rotation system, we must recover the number of faces of $G \rightarrow S$, hereandafter denoted $|F(G\rightarrow S)|$, by using a face-tracing algorithm as described more fully in \cite[\S 3.2.6]{GT}; the cyclic orderings of the incidence of edges at vertices may be used to reconstruct the facial boundary walks themselves, and account for each of the two times each edge must appear in the boundary of a $2$-complex homeomorphic to a surface.  (The reader is encouraged to consider Figure \ref{fig:rotationSystems} in light of this description.)  Since the contraction of an edge of $G\rightarrow S$ results in a cellular embedding of a different graph in $S$ with the same number of faces as $G\rightarrow S$, one may accelerate the determination of the number of faces of $G\rightarrow S$ by contracting a spanning tree $T$ of $G$;  upon the contraction of each edge link $e$ of $T$, the rotations on each of the end vertices of $e$ are spliced together, and if $e$ is a type-$1$ edge, then one of the spliced rotations is reversed to account for the untwisting of $e$ that happens while $e$ is contracted.  With $G\rightarrow S$ reduced to a one-vertex graph embedding, we now have a single finite list to analyze, and the face-tracing algorithm will go more quickly.  With the requisite information in hand, given a roation system of $G\rightarrow S$, we may determine the number of handles and crosscaps of $S$ by noting that a surface $S$ with $h$ handles, $c$ crosscaps, and Euler characteristic $\chi(S)$ obeys the equation \begin{equation}\label{equation:eulerchar} \chi(S) =|V|-|E|+|F(G\rightarrow S)| = 2-2h-c.\end{equation}

The piecewise linear encoding of a $2$-manifold is now given as a rotation system, which is in terms of tuples that are locations of vertices and the cyclic orderings of signed edges incident to them.

\begin{remark}The reader who is familiar with persistent homology may expect that we will be treating triangulations encoded by rotation systems.  This is not always going to be the case.  While the techniques of persistent homology do require that the $2$-cells of cellular complexes be triangles, a rotation system does not.  A rotation system may have faces bounded by any number of edges, and a face may contain both sides of an edge in its boundary, as in Figure \ref{fig:rotationSystems}.\end{remark}

\section{Some necessary results from manifold learning}\label{section:manifoldBackground}

As mentioned above, our algorithm for constructing the $2$-complex requires estimating the tangent plane of the embedded surface.  In order to take a rigorous approach, we will employ a method \cite{Berry} which provably estimates the tangent directions of the manifold $M$ at a point $x\in f(M)$ in the limit of large data.  Intuitively speaking, a data set $\{x_i\}_{i=1}^N \subset f(M)$ which lies on the embedded manifold will be approximately linear in a local region (meaning some sufficiently small $\epsilon$-ball) around any given point, $x$.  We can find a basis for this linear space using a singular value decomposition of the vectors pointing from $x$ to the nearby data points $x_i$ and choosing the singular vectors with the largest singular values.  Of course, this will only be true given a sufficiently large data set that the nearest neighbors of each point are close enough to the point $x$ so that the difference $x_i-x$ is close to being a tangent vector.  Moreover, there is a question of `how large is large enough' for the singular values.  These issues are addressed using the method introduced in \cite{Berry} which derives a provable cutoff based on weighting the vectors $x_i-x$ and using weights that decay exponentially in the length of the vector.  

\begin{figure}[h]
\begin{center}
\includegraphics[width=.9\linewidth]{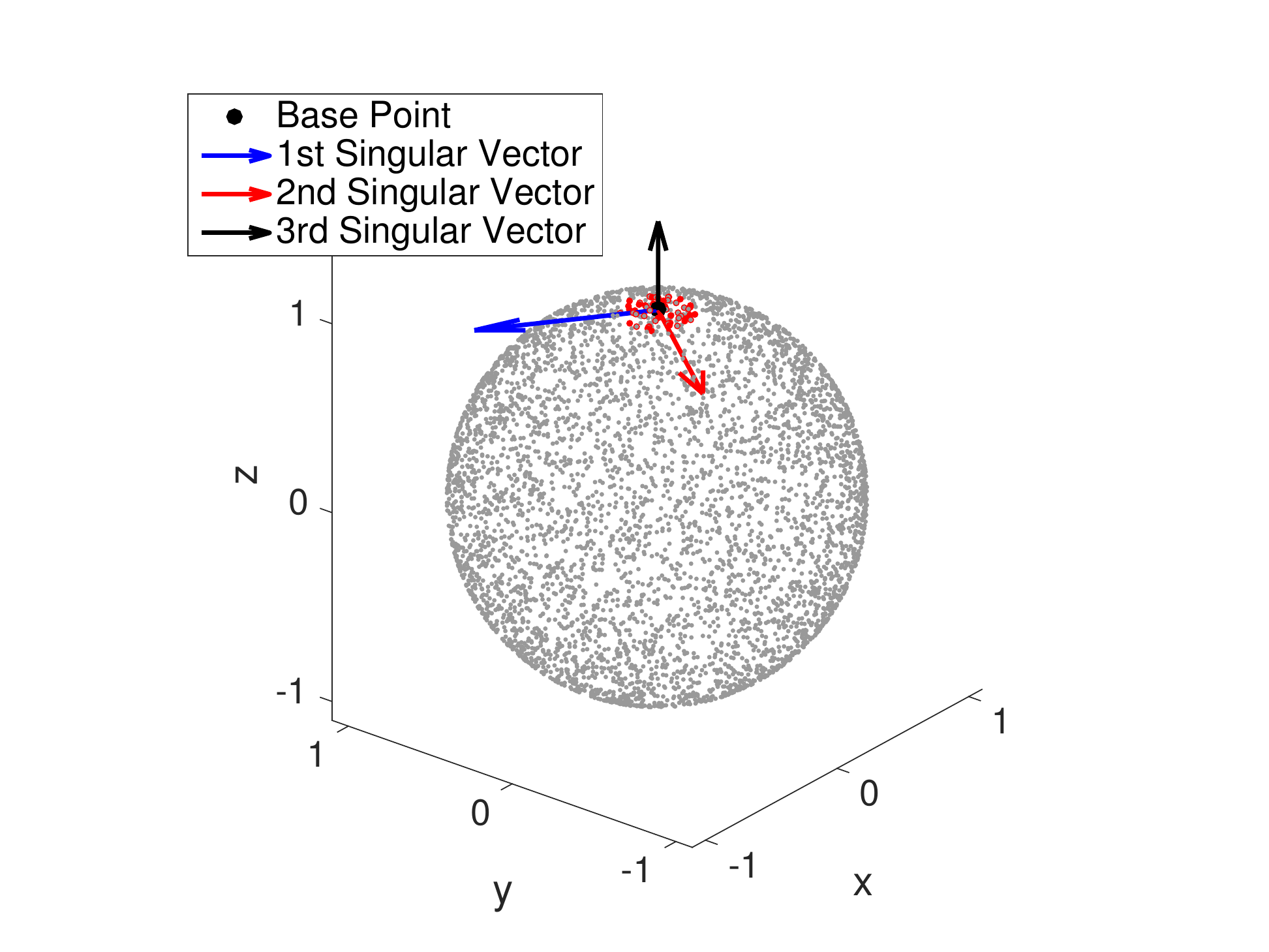} \\
\includegraphics[width=.45\linewidth]{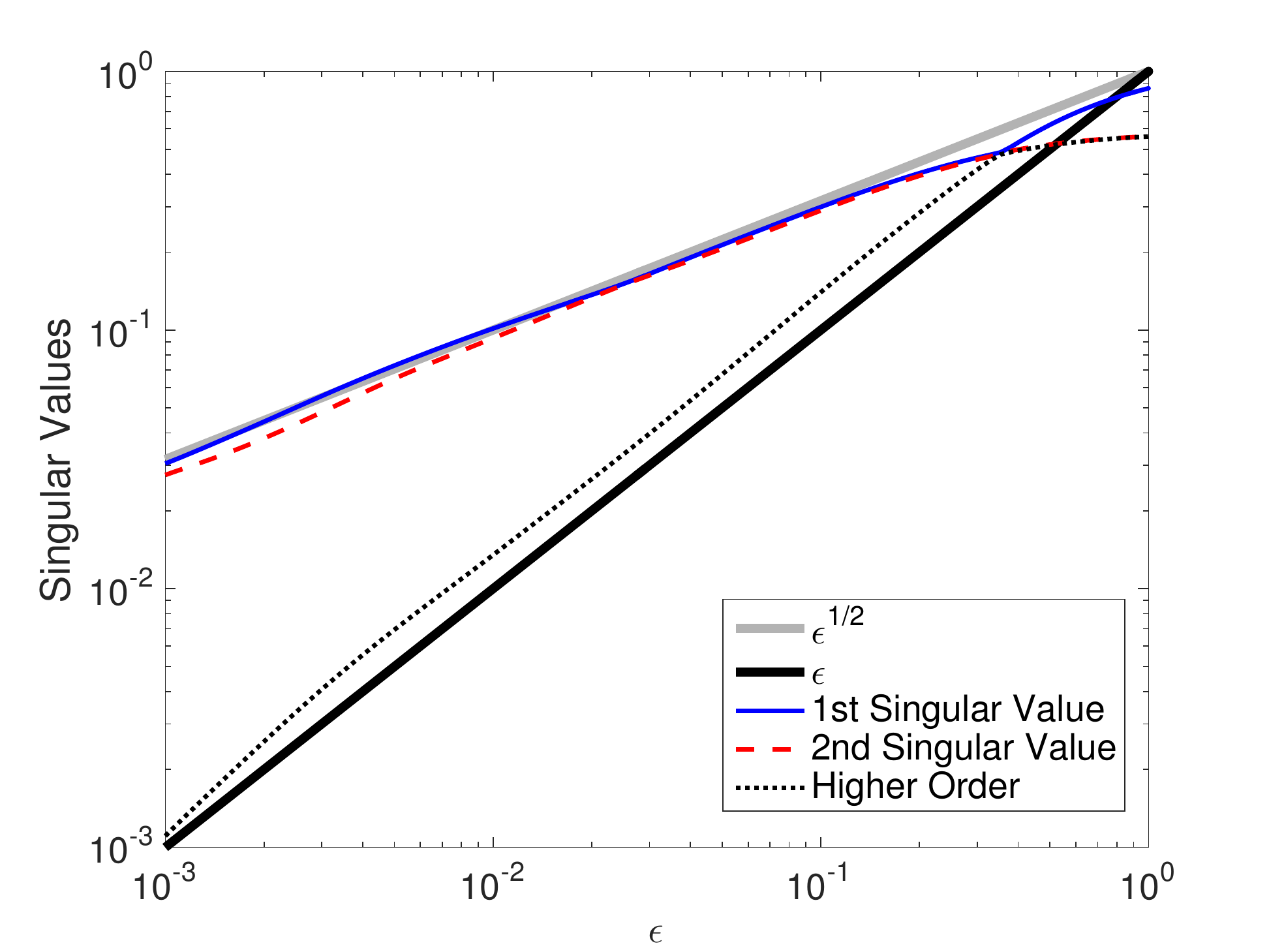}
\includegraphics[width=.45\linewidth]{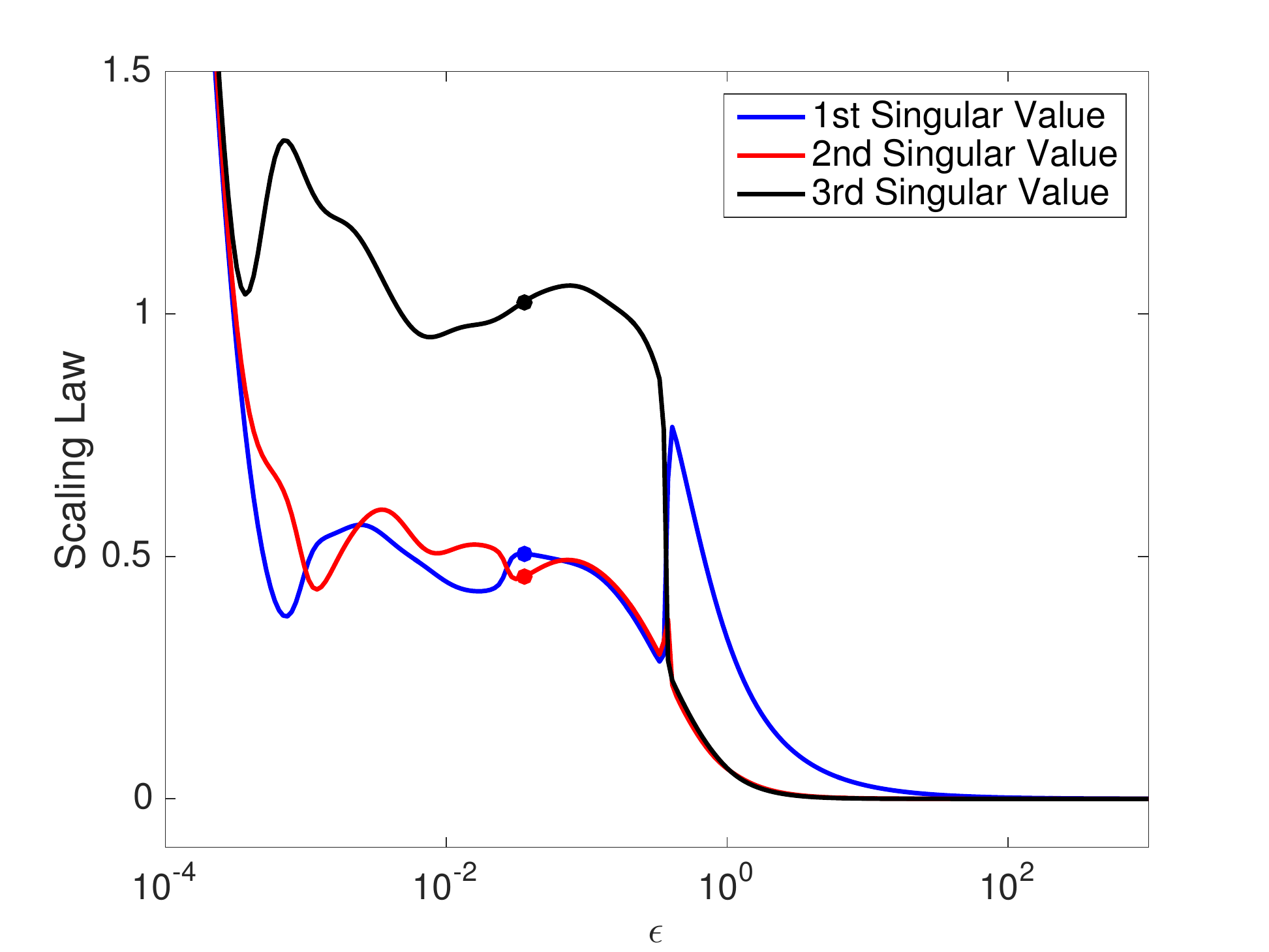}
\caption{\label{fig:scalinglaws} Top: A data set on the sphere in $\mathbb{R}^3$ showing a chosen base point and the singular vectors of the matrix $X_{\epsilon}$ (with $\epsilon$ as shown in the right figure) the red data points are inside the $\epsilon$ neighborhood.  Bottom, left: the singular values of the matrix $X_{\epsilon}$ as a function of $\epsilon$ showing the $\sqrt{\epsilon}$ and $\epsilon$ scaling laws of the singular values corresponding to the tangent and normal vectors respectively.  Bottom, right: extracting the scaling laws $\alpha_j$ for each singular value as a function of $\epsilon$.  The value $\epsilon$ is automatically chosen so that all the scaling laws are simultaneously as close to stationary as possible (see \cite{Berry} for details).}
\end{center}
\end{figure}

We assume that $\{x_i\}_{i=1}^N$ are independent samples from a distribution with a smooth density $q(x)$ defined on $f(M)$.  Fix a point $x \in f(M)$ and for simplicity we will weight all of the vectors pointing from $x$ to a data point $x_i$, so we define the vectors
\[ w_i(\epsilon) = \exp\left(-\frac{||x_i-x||^2}{4\epsilon}\right)(x_i - x) \]
(where $\epsilon$ is the neighborhood radius that will be fixed later but is allowed to very in this section).   We also define the normalization factor 
\[ D(x,\epsilon) = \sum_{i=1}^N \exp\left(-\frac{||x_i-x||^2}{2\epsilon}\right) \]
and we let $X_{\epsilon}$ be the matrix whose $i$-th column is the weighted vector $D(x,\epsilon)^{-1/2}w_i(\epsilon)$.  Then Corollary 3.2 of \cite{Berry} shows that
\[ \lim_{N\to\infty} X_{\epsilon}X_{\epsilon}^\top = \epsilon Df(p)^\top Df(p) + \mathcal{O}(\epsilon^2) \]
where $x = f(p)$.  The dependence of $X_{\epsilon}X_{\epsilon}^\top$ on $N$ is due to $X_{\epsilon}$ having $N$ columns.  Notice that the embedding $f$ maps the manifold into a higher dimensional Euclidean space, and so the range of $Df(p)$ is the image of the tangent space $T_pM$ in the embedding space.  So, both the row and column space of $Df(p)^\top Df(p)$ span the image of the tangent space, and the components orthogonal to the tangent space will be contained in the higher order term $\mathcal{O}(\epsilon^2)$.  In other words, given a singular vector $v$ (with $||v||=1$) of the matrix $X_{\epsilon}$, the singular value is given by 
\[  \lim_{N\to \infty} \sigma_v(\epsilon) = \lim_{N\to \infty} \sqrt{v^\top X_{\epsilon}X_{\epsilon}^\top v} = \sqrt{\epsilon v^\top Df(p)^\top Df(p) v} + \mathcal{O}(\epsilon).  \]
So if $v$ lies in the tangent space (namely the image of $Df(p)$), then in the limit of large $N$ and small $\epsilon$ the singular value will be order-$\sqrt{\epsilon}$, whereas if $v$ is orthogonal to the tangent space then the singular value will be order-$\epsilon$.  In order to differentiate these two cases, we construct $X_{\epsilon}$ and compute the singular value decomposition for a range of values of $\epsilon$, so that the singular vectors $v_j(\epsilon)$ and singular values $\sigma_j(\epsilon)$ are functions of $\epsilon$ (where $j=1,...,n$ the dimension of the embedding space).  We then estimate whether each singular value scales like $\sqrt{\epsilon}$ or $\epsilon$ by estimating the derivative
\[ \alpha_j(\epsilon) = \frac{d\log(\sigma_j(\epsilon))}{d\log\epsilon}. \]
Finally, we sort the vectors according to whether $\alpha_j(\epsilon) \approx 1/2$ (which are the tangent vectors) or $\alpha_j(\epsilon) \geq 1$ (which are the normal vectors).  This process can be applied to find a good approximation of the tangent space of an embedded manifold of any dimension given sufficient data.

In Figure \ref{fig:scalinglaws} we demonstrate this process for a data set sampled from a uniform distribution on the unit sphere in $\mathbb{R}^3$.  The data set was produced by sampling $5000$ random points in $\mathbb{R}^3$ and projecting them onto the sphere by dividing each point by its 2-norm.  For clarity, we chose the base point to be the point closest to the North pole, and then we applied the method described above for a wide range of $\epsilon$ values, producing the singular values of the matrix $X_{\epsilon}$ for each $\epsilon$.  In the bottom left figure we show how the singular values are either proportional to $\sqrt{\epsilon}$ or $\epsilon$.  As shown in the top figure, the $\sqrt{\epsilon}$ singular values correspond to vectors lying in the tangent space while the $\epsilon$ singular values are normal to the tangent space.  Finally, to show how the selection process can be automated, we plot the scaling laws $\alpha_j$ for each singular value in the bottom right panel of Figure \ref{fig:scalinglaws}.

Of course, to determine the tangent vectors we must choose a particular value for $\epsilon$.  A method of choosing an $\epsilon$ value is given in \cite{Berry}, but in our implementation below we used a faster alternative which simply takes $\epsilon$ to be the distance to the $k$-th nearest neighbor (in this case we simply choose $k=20$ but more generally $k$ should typically increase as the number of data points increases). From this point on $\epsilon$ be fixed for each base point $x$.

\section{The Algorithm}\label{section:theAlgorithm}

At the highest level our approach decomposes into three steps:
\begin{enumerate}
\item{construct a rotation system encoding a $2$-complex from the data set represented by a directed graph, and}
\item{contract a spanning tree, and} 
\item{compute the Euler characteristic of the resulting $2$-complex.}
\end{enumerate} 
We will cover these three parts of the algorithm in three separate subsections since each part of the algorithm may be of independent interest.

\subsection{Constructing the 2-complex}\label{section:construction}

\begin{figure}[h]
\begin{center}
\includegraphics[width=.49\linewidth]{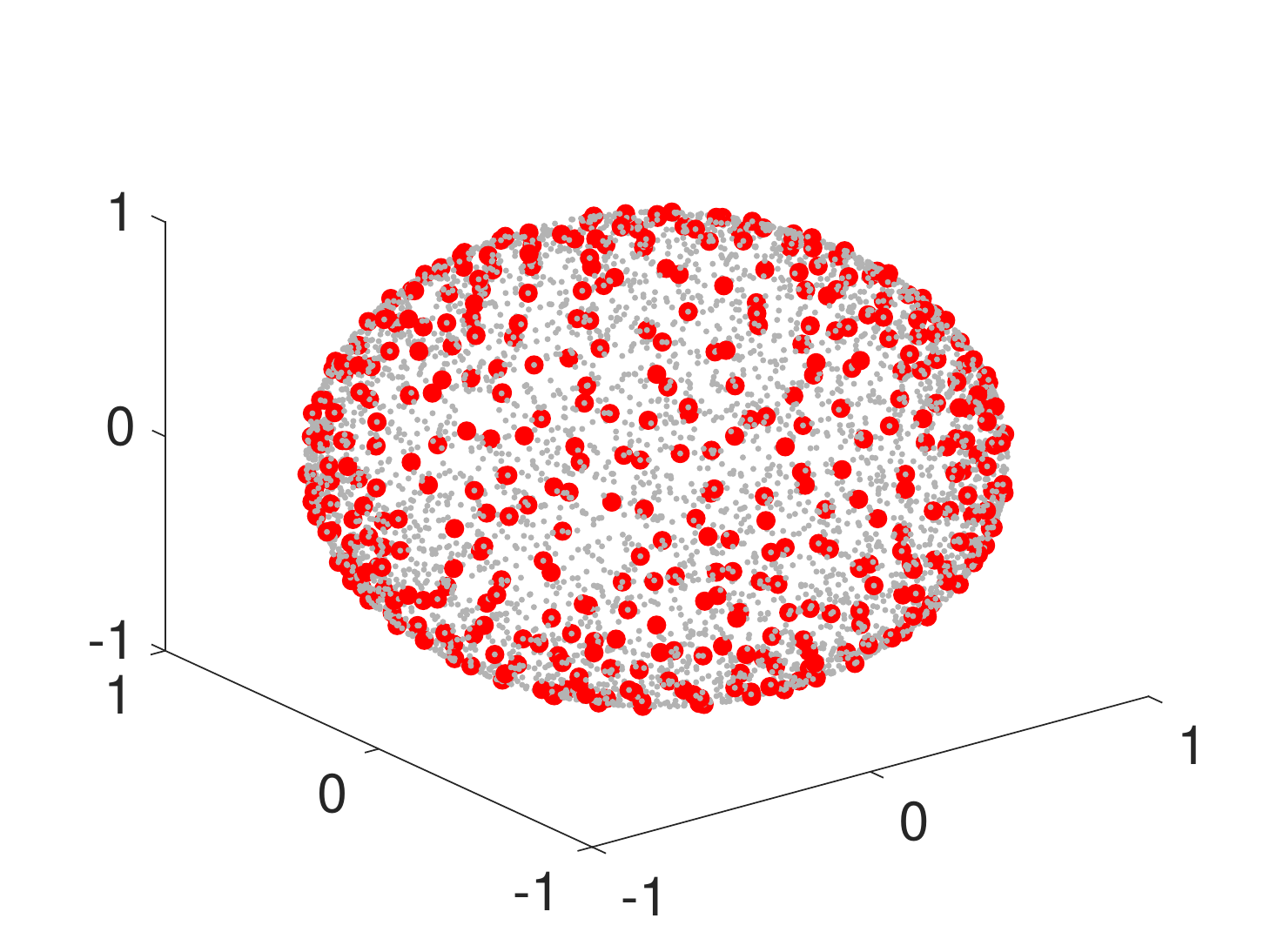}
\includegraphics[width=.49\linewidth]{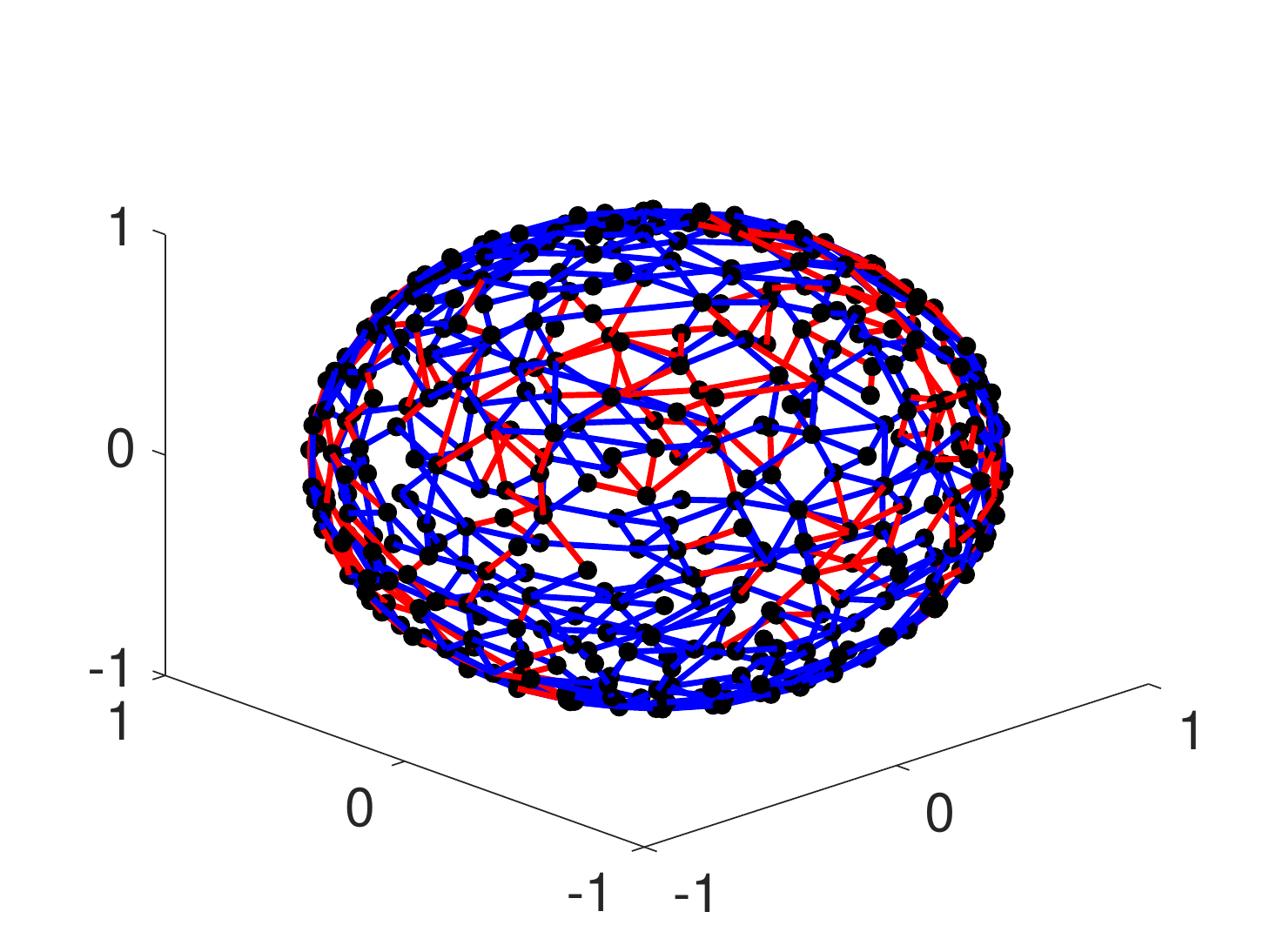}
\caption{\label{fig:triangulation} Left: The full $5000$ data points (small grey points) along with the vertices $V$ selected by the algorithm (large red points) on the sphere in $\mathbb{R}^3$.  Right: The cell complex on the sphere produced by the algorithm, where color indicates twisted edges (red) meaning the adjacent vertices have opposite orientations on their respective tangent space bases, or untwisted edges (blue) meaning the orientations are the same.  Since the sphere is orientable, any path from a point back to itself should have a net zero twist.}
\end{center}
\end{figure}

The basic strategy for constructing the 2-complex is to start from any point on the manifold and to use the tangent vectors to traverse the manifold, selecting a subset of the data points as vertices in a graph and adding as many edges as possible as we go.  We should note that by representing our 2-complex using only a subset of the data points as vertices, there is a significant reduction in the complexity of the data set.  This is illustrated in Fig.~\ref{fig:triangulation} where we show the full data set and the subset of vertices chosen to represent the data, as well as the final cell complex.  Each vertex in the 2-complex will be `summarizing' the neighboring points in the data set, meaning that not every data point will need to be a vertex in the cell complex, but each data point will be within a small neighborhood of a vertex point.  The algorithm will terminate when every point from the data set is either a vertex of the resulting 2-complex or is discarded for being sufficiently close to a vertex.  The most challenging part of this algorithm is determining when we can validly add an edge.  The central challenge is that in order to respect the 2-complex structure the edges cannot be allowed to cross.  If two edges cross, the the result of this is that the rotation system may falsely interpret the result as indicating the presence of an additional handle on the surface, thereby incorrectly increasing the Betti number.  This runs contrary to the idea of highly local approximations giving rise to a correct global approximation of a surface.  In order to test for crossings we will project the edges into the tangent space and apply a geometric test.  We now detail the steps of the algorithm.

Our algorithm builds the following data structures.  
\begin{itemize}
\item A list of vertices $V$ containing the indices of the subset of data points which are the vertices in our complex.  
\item A list $C$ containing one integer for each data point which gives the index of the vertex which summarizes that data point.  
\item A collection of $2\times n$ matrices (one matrix for each entry in $V$) containing the coordinates of two orthogonal vectors in $\mathbb{R}^n$ that span the tangent space at the corresponding vertex from $V$.
\item  An adjacency matrix $A$ which for each edge in the complex contains the angle (between $[0,2\pi)$) of the vertex with respect to the two orthogonal vectors at the base vertex (so $A$ is non-symmetric).  The sign of each entry in $A$ indicates whether the two tangent spaces have the same (positive) or opposite (negative) orientations.  
\end{itemize}

We start by choosing an arbitrary initial point from the data set and adding this point as a vertex in our graph by adding its index to the list $V$.  We then apply the method from Section \ref{section:manifoldBackground} to estimate the tangent vectors to the surface at this point.  This method also returns a radius $\epsilon$ inside which the tangent space is deemed close to the manifold (in the sense that the weighted neighbors are well approximated by a plane).  We mark all of the data points within this radius as `summarized' by writing the index of the vertex into the entries of $C$ corresponding to these data points.  To expand the graph we then follow each of the tangent vectors in the positive and negative directions.  In each direction we find the nearest data point outside of the $\epsilon$ radius; if that data point is not summarized then we recursively call the above procedure (adding that data point to the vertices, finding the tangent space, and exploring the tangent directions).  If the data point is already summarized by a vertex, then we need to determine if we can add an edge between the two vertices.  The issue is that this edge could cross an existing edge (meaning that the corresponding line segments may geometrically cross when projected into the tangent space of one of the vertices). 

To check for crossings, we first project the potential edge into the tangent space of the base vertex, and then using the adjacency matrix $A$ we consider all of the edges that connect to vertices within $10\epsilon$ of the two vertices which would be connected to the potential edge.  Each of these nearby edges is also projected into the tangent space and we determine geometrically if the two edges cross.  If the potential edge does not cross any nearby edges, then we compute the angles of the edge relative to the tangent vectors of both the adjacent vertices.  Finally, to avoid small angles in our complex, we check that the angle between the potential edge and all the existing edges on each of the adjacent vertices are greater than a specified minimum ($\pi/6$ in our implementation); this will be discussed further in Section \ref{section:smallAngles}.  If all of these criteria are met, we can add the edge to our graph by inserting the computed angles into the adjacency matrix.  We should note that the adjacency matrix is non-symmetric since the each edge will have two associated angles (one angle for each adjacent vertex, representing the angle of the edge projected into the chosen tangent space coordinates).  

Finally, since the orientations of the tangent spaces may be different, we compute the sign of the determinant of the linear map that takes the tangent vectors at one of the adjacent vertices to the tangent vectors at the other adjacent vertex, and we use this determinant as the sign of the two new entries in the adjacency matrix.  In other words, if the entry in the adjacency matrix is negative, that denotes that the tangent space coordinates of the adjacent vertices have opposite orientations (if the entry is positive they have the same orientation).

The above procedure continues recursively until no more edges or vertices can be recursively added.  At this point we check to see if every data point has been summarized, and if not, we randomly select an unsummarized data point to add and repeat the above procedure on this new point.  We continue this process until every data point has been summarized.  

\begin{figure}[h]
\begin{center}
\includegraphics[width=.9\linewidth]{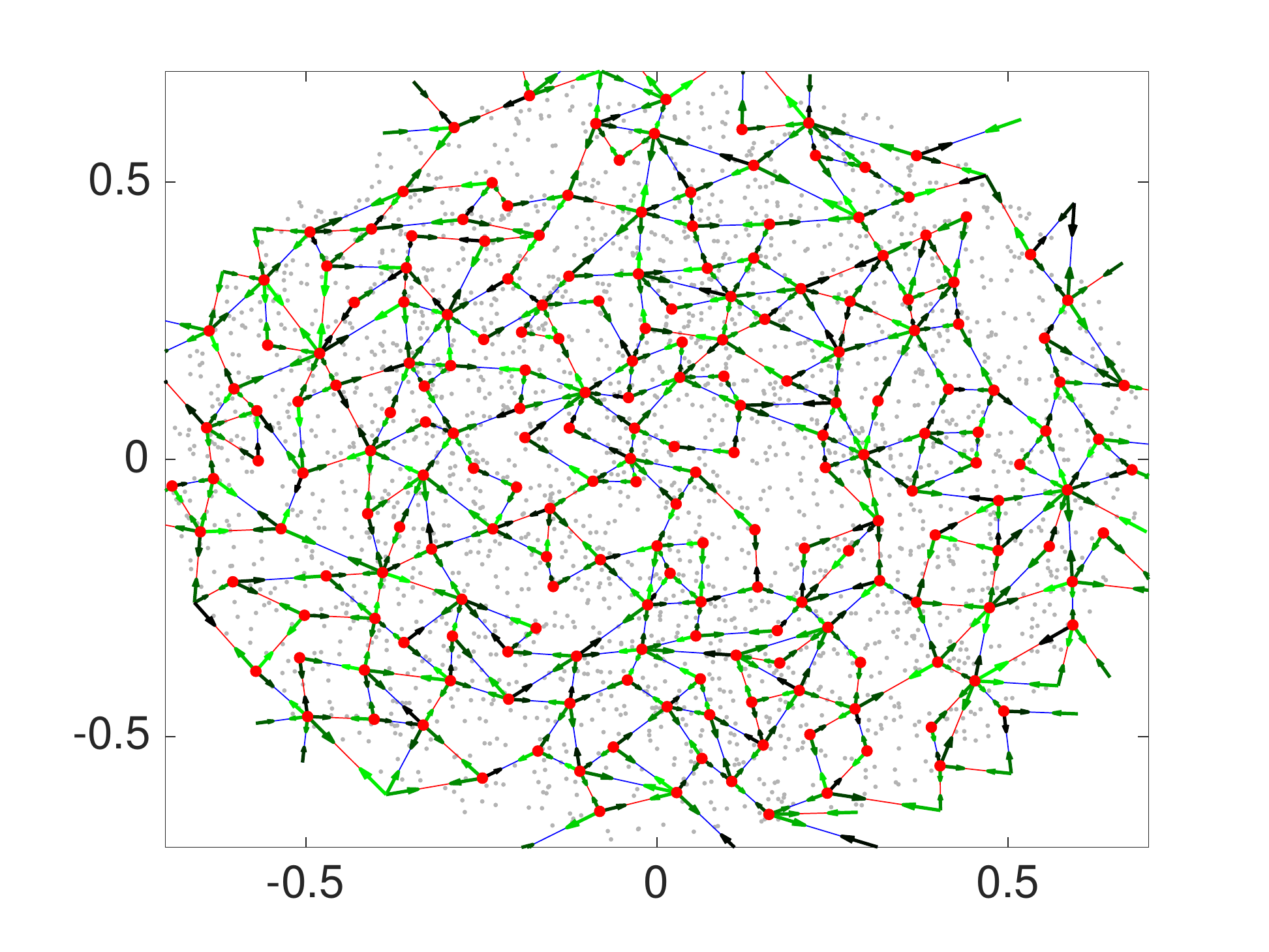}
\caption{\label{fig:rotationsystem} The rotation system is shown on an unfolded hemisphere of the sphere.  The original data points are shown as small black points and the vertices $V$ are shown as large red points.  Arrows are not shown to indicate directed edges, rather they are there to indicate the angle relative to the tangent space of the different vertices, which in turn determines the rotation system.  As the color changes from lighter to darker the angle increases in counter-clockwise order.  A thinner edge is shown connecting every pair of adjacent vertices (visible only between the arrows) which is colored either red for twisted edges or blue for untwisted edges.}
\end{center}
\end{figure}

An example of the final result of the algorithm is shown in Fig.~\ref{fig:rotationsystem} for the data set on the sphere.  For clarity we show only a hemisphere that has been unfolded into $\mathbb{R}^2$.  The vertices selected by the algorithm are shown as large red dots and the algorithm guarantees that each data point (small black dot) is within the $\epsilon$ neighborhood of a vertex.  The information in the adjacency matrix is shown by the green arrows, whose darkness indicates the angle of the edge relative to the coordinates defined by the tangent vectors based at the base vertex.  Notice that as one moves around a vertex, the darkness of the arrows has only a single discontinuity (indicating the jump from $2\pi$ back to $0$ radians).  Beneath each pair of arrows, the edge is also drawn as a thinner line which is colored red for twisted and blue for untwisted edges.  Notice that the direction of the darkness gradient of the arrows is reversed when comparing two vertices connected by a twisted (red) edge.  Finally, notice that all of the angles between edges are at least $\pi/6$ and there are no edges that cross one another.  These are two conditions which are crucial in the construction of the complex in order to insure a proper cell complex and associated rotation system.

\subsection{Building the rotation system}

At this point we have constructed a graph embedding into the manifold described by the data along with the supplementary information of the edge angles and whether each edge is orientation preserving or orientation reversing.  We can now easily sort the edge angles to produce a rotation system.  However, for ease of computation in the next steps, it is convenient to build a rotation system on a single vertex by collapsing a spanning tree.  Contracting the spanning tree also requires us to use local sign switches to push all of the orientation reversing edges outside of the spanning tree.  Using the data structures constructed in the previous step of the algorithm this can easily be accomplished.

Essentially we will be applying a depth-first approach to building a spanning tree, while carefully pushing orientation-reversing edges outside of the tree as we build it.  The rotation on each vertex will simply be an ordered list of numbers where each number corresponds to edge, and each number occurs exactly twice, once for each end of the corresponding edge (as illustrated in Fig. \ref{fig:rotationSystems}).  The ordering of the list indicates the ordering of these edges in the single vertex graph, which remains after the spanning tree is contracted.  In order to create a number corresponding to each edge that is independent of the vertex that is the base, we simply multiply the smaller of the vertex numbers by the number of vertices and then add the larger of the two vertex numbers.  The edge number will either be positive or negative depending on whether the edge is orientation preserving or reversing.  We construct the rotation system as follows in the next paragraph. 

Starting from any vertex of the graph mark the vertex as visited and then, we sort the edges according to their angles (which are stored in the adjacency matrix) and iterate through each edge in order.  If an edge goes to a visited vertex, then we add the edge number to the rotation system, including the sign that indicates the edge is orientation preserving or reversing.  If an edge goes to an unvisited vertex then we `visit' that vertex recursively.  In order to visit a vertex, we first check if the edge is orientation reversing, and if so then we perform a local sign switch on that vertex.  Next, we shift the numbering of the edges so that the edge connecting the previous and current vertices is the first edge in the ordering, and we then remove that edge.  We can now proceed recursively, first marking this vertex as visited and then traversing each edge in the new ordering, either adding the edge to the rotation system if the other adjacent vertex is visited or recursively `visiting' the vertex.  This procedure will terminate automatically when every vertex has been visited and the output will the be the rotation system (ordered list of signed edge numbers).

\subsection{Computing Euler characteristic and classifying the resulting surface}

\begin{figure}[h]
\begin{center}
\includegraphics[scale=0.5]{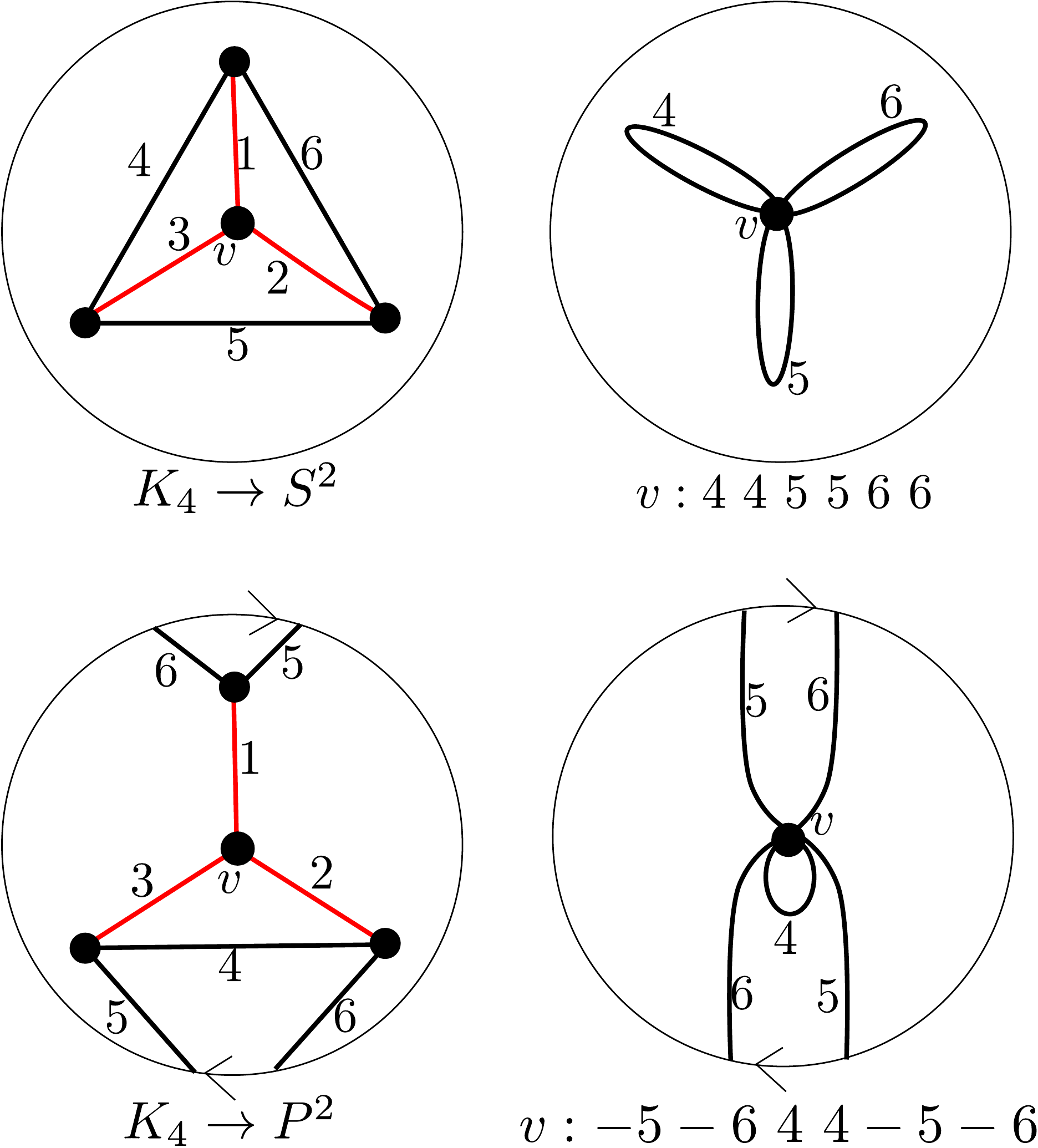}
\caption{Two examples of the contraction of a spanning tree of cellular graph embeddings, and the resulting rotation systems.  Each list is a cyclic ordering of edge ends incident to the vertex.}\label{fig:ContractedTree}
\end{center}
\end{figure}

At this point we have constructed a graph embedding and represented it in the form of a rotation system.  As mentioned previously and illustrated in Figure \ref{fig:ContractedTree}, we may contract a spanning tree and produce a rotation system of an embedding of a one-vertex graph in the same surface.  Since the Euler characteristic is invariant under homotopy, we may evaluate the Euler characteristic of the $2$-complex created by the contraction.  Note that here, because the result of contracting a spanning tree is a cellular embedding of a bouquet of loops, we will denote a cyclic ordering of edge ends incident to the vertex. Positively-signed integers indicate orientation-preserving edges, and negatively-signed integers indicate orientation-reversing edges.

To evaluate the Euler characteristic of a surface encoded in the form of a rotation system, one may utilize the face tracing algorithm described in Section \ref{section:graphBackground} and in \cite[\S 3.2.6]{GT}.  The reader is encouraged to note that the benefit of contracting edges is now apparent in that one needs to search through only one list of integers instead of several such lists in order to find the number of faces of the encoded embedding.

With the Euler characteristic in hand, all one must do to classify the surface is to apply Equation \ref{equation:eulerchar} with the orientability of the surface obtained by checking if there are any negatively signed integers in the rotation system; moreover, if the Euler characteristic is odd, then we already know that the surface is non-orientable.  For example, if $\chi(S)=0$ and $S$ is known to be non-orientable, then $S$ must be a surface with two crosscaps (the Klein bottle) and not the torus.  If $\chi(S)=-1$, then Equation \ref{equation:eulerchar} makes it clear that since the addition of a handle to a surface lowers the Euler characteristic by $2$, then $S$ must be the three-crosscap (Dyck's) surface.

\begin{figure}[h]
\begin{center}
\includegraphics[width=.45\linewidth]{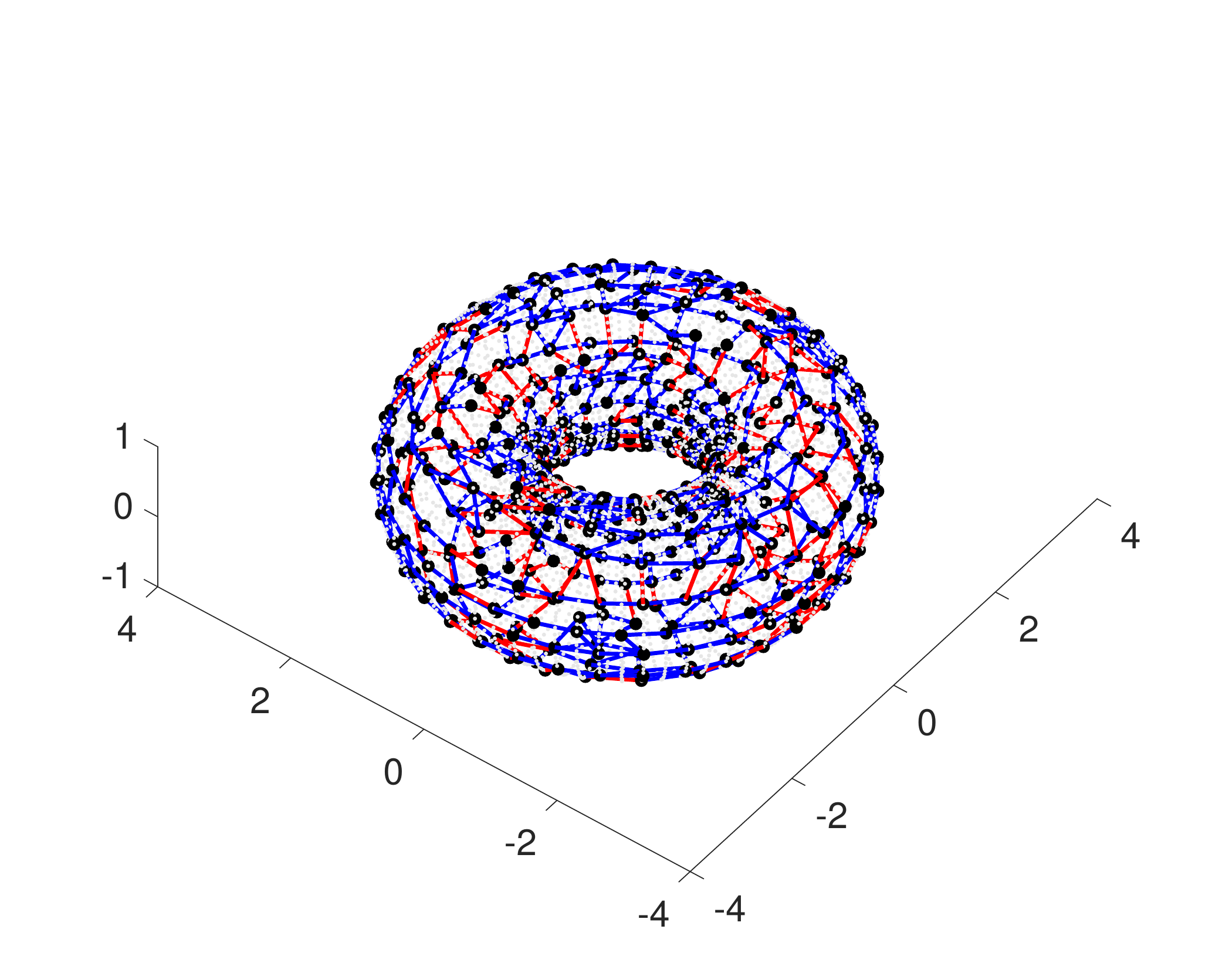}
\includegraphics[width=.45\linewidth]{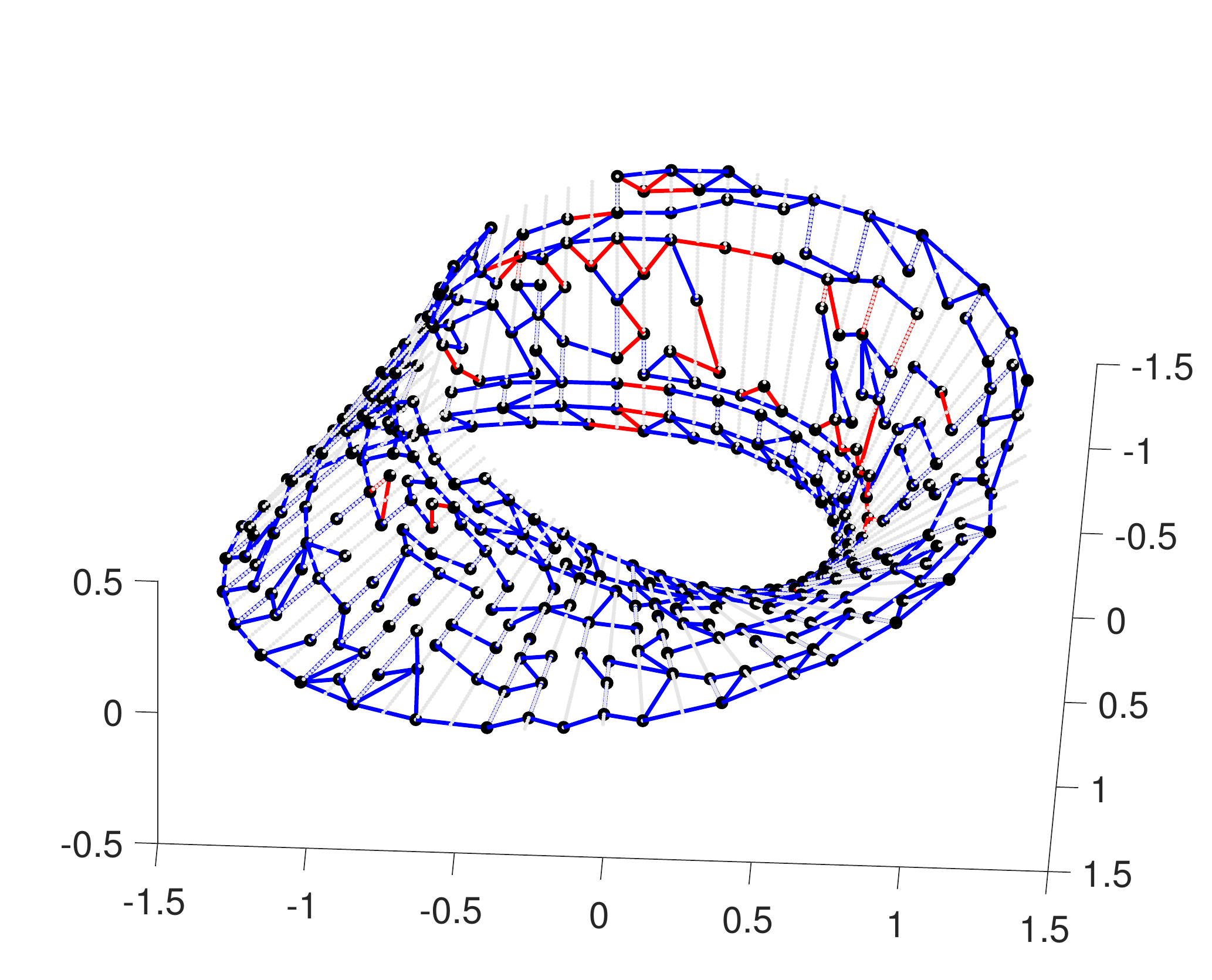}\\
\includegraphics[width=.45\linewidth]{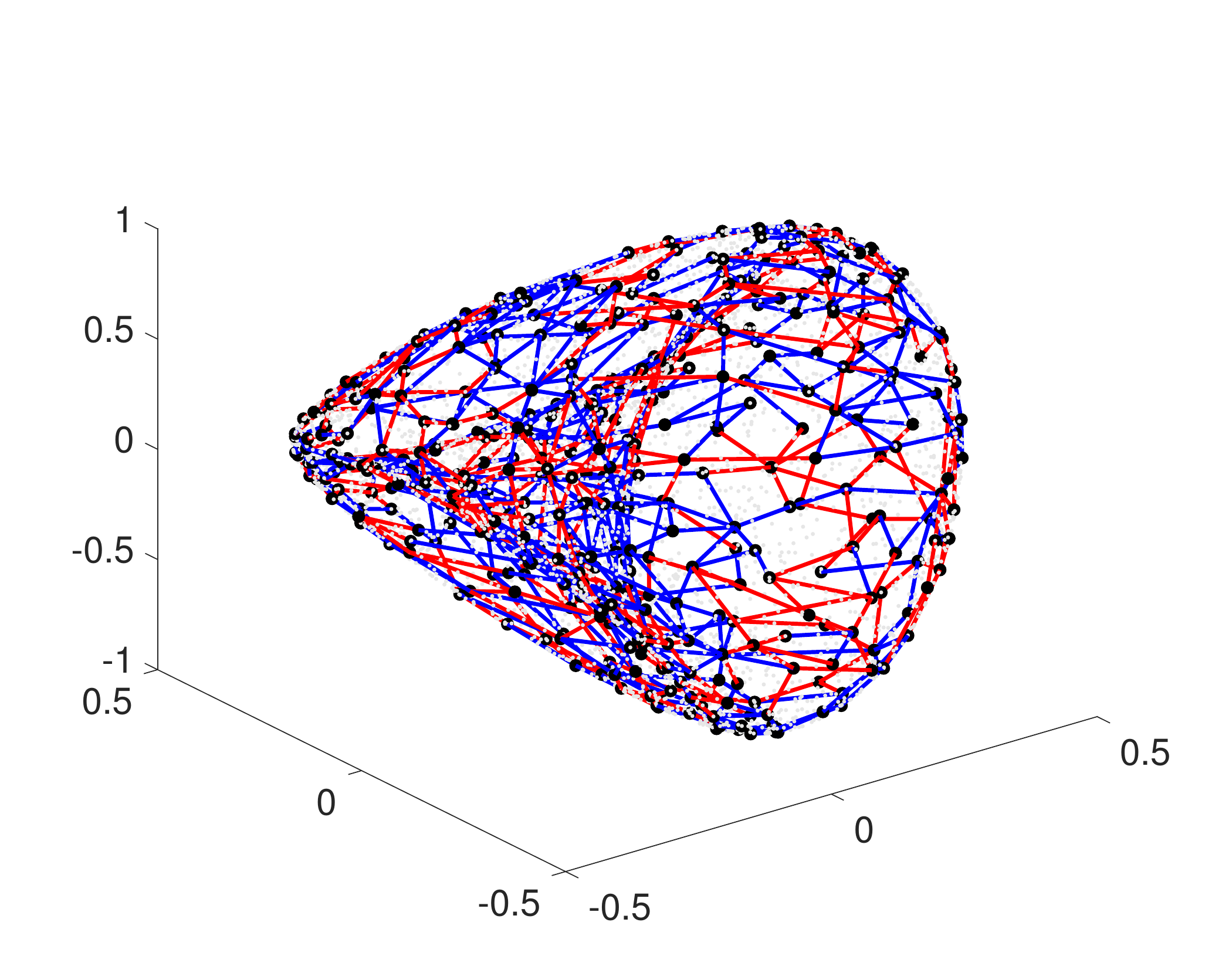}
\includegraphics[width=.45\linewidth]{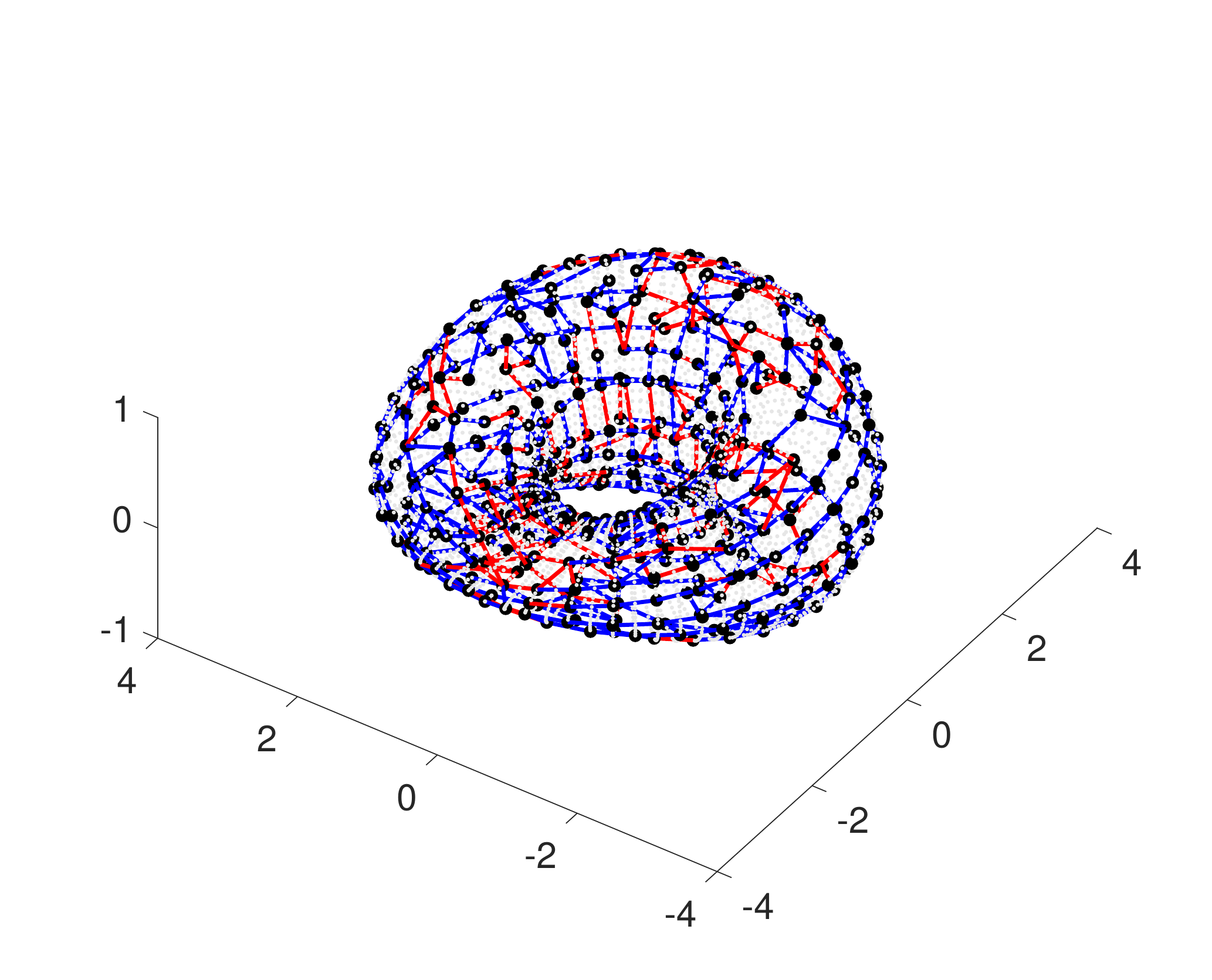}\\
\includegraphics[width=.9\linewidth]{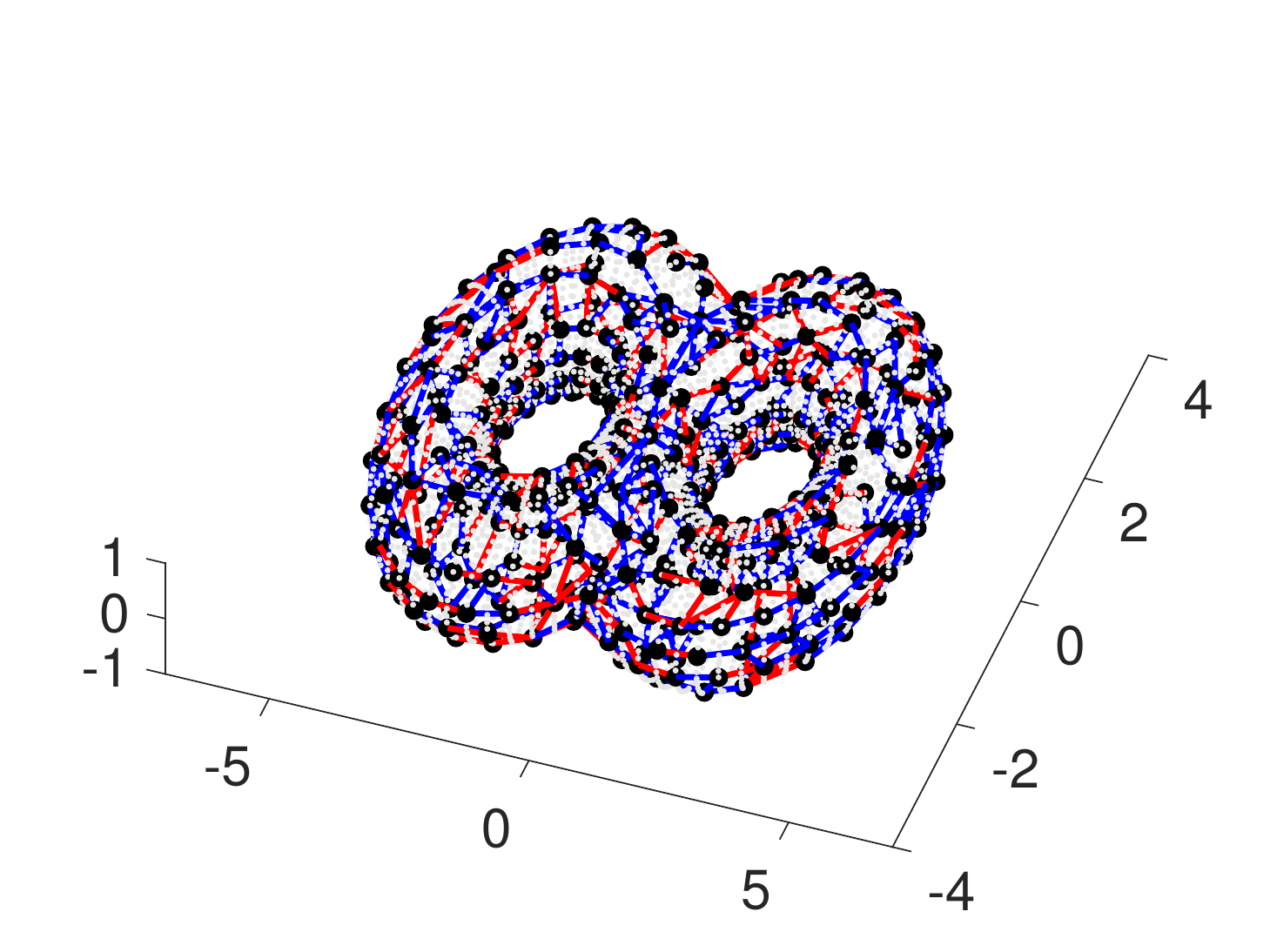}
\caption{\label{fig:eulercharacteristics} In addition to the sphere, our method was reliably able to obtain the correct Euler characteristic on the above data sets representing a torus, Mobius band (seen as a projective plane by our algorithm), RP2, Klein bottle, and a genus-2 torus.  For each data set we show the cell complex produced by our algorithm (projected into $\mathbb{R}^3$).}
\end{center}
\end{figure}

\section{Experimental data and closing remarks}\label{section:data}

\subsection{Experimental Data}\label{section:experimentalData}

We ran our Matlab code on $100$ noisily sampled point clouds from embeddings of each of the following surfaces, one of which has a boundary component (which we will explain later in Section \ref{section:boundary}): the sphere, the projective plane, the torus, the M\"obius band, the Klein bottle, and the genus-2 torus.  Our method produced encouraging results, which are delineated in Figure \ref{figure:dataTable}.

\begin{figure}[h!]

\begin{tabular}{| c | c | c | c | c | c | c | c | c | c | c |}
\hline

\multicolumn{1}{|c}{ } & \multicolumn{2}{r|}{Orientable} & \multicolumn{8}{c|}{Euler Characteristic }\\ 
\hline 
Surface & Yes & No & -6 & -4 & -3 & -2 & -1 & 0 & -1 & 2\\
\hline
sphere & 100 & 0 &  0&    0  &   0 &    3    & 0  &  21   &  0  &  76\\
Torus & 100  &   0 &  1  &  1 &    0   & 10 &    0  &  54  &   0  &  35\\
M\"obius band & 1   & 99 &  0 &  0  &   1  &   0 &    8  &   0  &  90 &    1 \\
RP2 & 0  & 100 & 0  & 0  &   8  &   0  &  32  &   0  &  60  &   0 \\
Klein bottle & 0  & 100  &  0 &    0  &   0  &   9 &    0  &  90  &   1  &   0 \\
genus-2 torus & 100  &   0 & 3 &   11 &    0  &  86   &  0  &   0   &  0  &   0 \\
\hline
\end{tabular}
\caption{A table indicating the results produced from 100 noisily sampled embeddings of the surfaces (one with boundary) named in the table.  The results corresponding to the row containing the word ``sphere" indicate that the sphere was identified as being orientable $100$ of $100$ times, and having the proper Euler characteristic $76$ of $100$ times.}\label{figure:dataTable}
\end{figure}

\subsection{Surfaces with boundary}\label{section:boundary}
If the reader is thinking that a rotation system must encode a surface without boundary, then the reader is correct.  Indeed, as one's eyes follow the cyclic ordering of edges incident to a vertex, one's eyes are tracing a circular path through a disc that encapsulates that vertex.  However, one may still ask the question:  were our methods to be working correctly, what would they produce if we began with a point cloud sampled from an embedded $2$-manifold with boundary?  The answer is that the resulting rotation system would fail to encode the boundary component, thus filling it in to complete a surface without boundary.  An encouraging result was that given point clouds sampled from M\"obius bands, our experimentation reliably produced rotation systems that encoded projective planes, as evidenced by Figure \ref{figure:dataTable}.

\subsection{On the importance of adding angles that are not too small}\label{section:smallAngles}
\begin{figure}[h!]
\begin{center}
\includegraphics[scale=0.5]{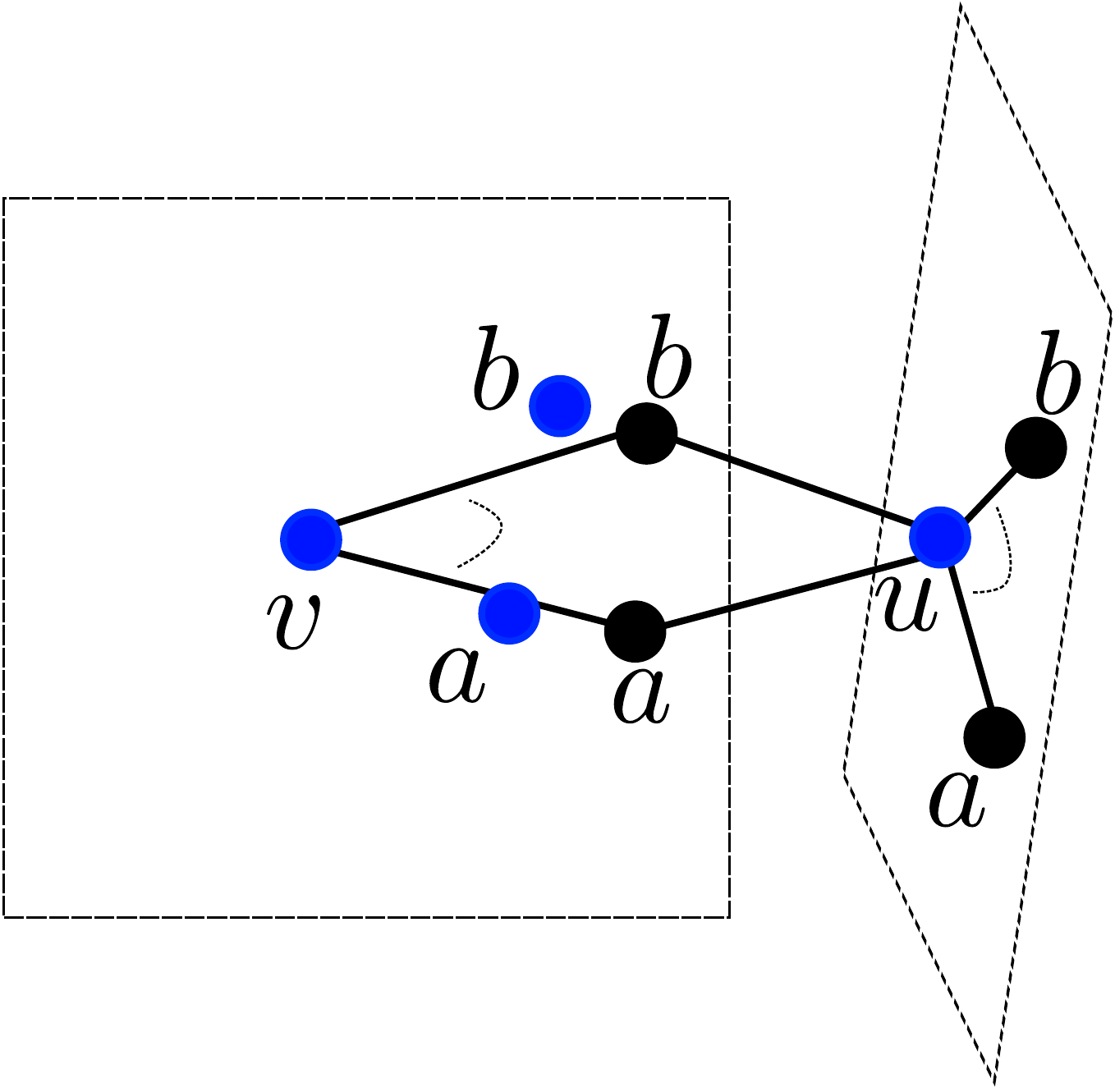}
\caption{An example of the inconsistent cyclic ordering of vertices seemingly caused by small corner angles and projections into tangent planes that intersect in angles of significance.  The projected vertices are black, and the nonprojected vertices are blue.}\label{figure:RotationsReversed}
\end{center}
\end{figure}

\begin{figure}[h!]
\begin{center}
\includegraphics[scale=0.5]{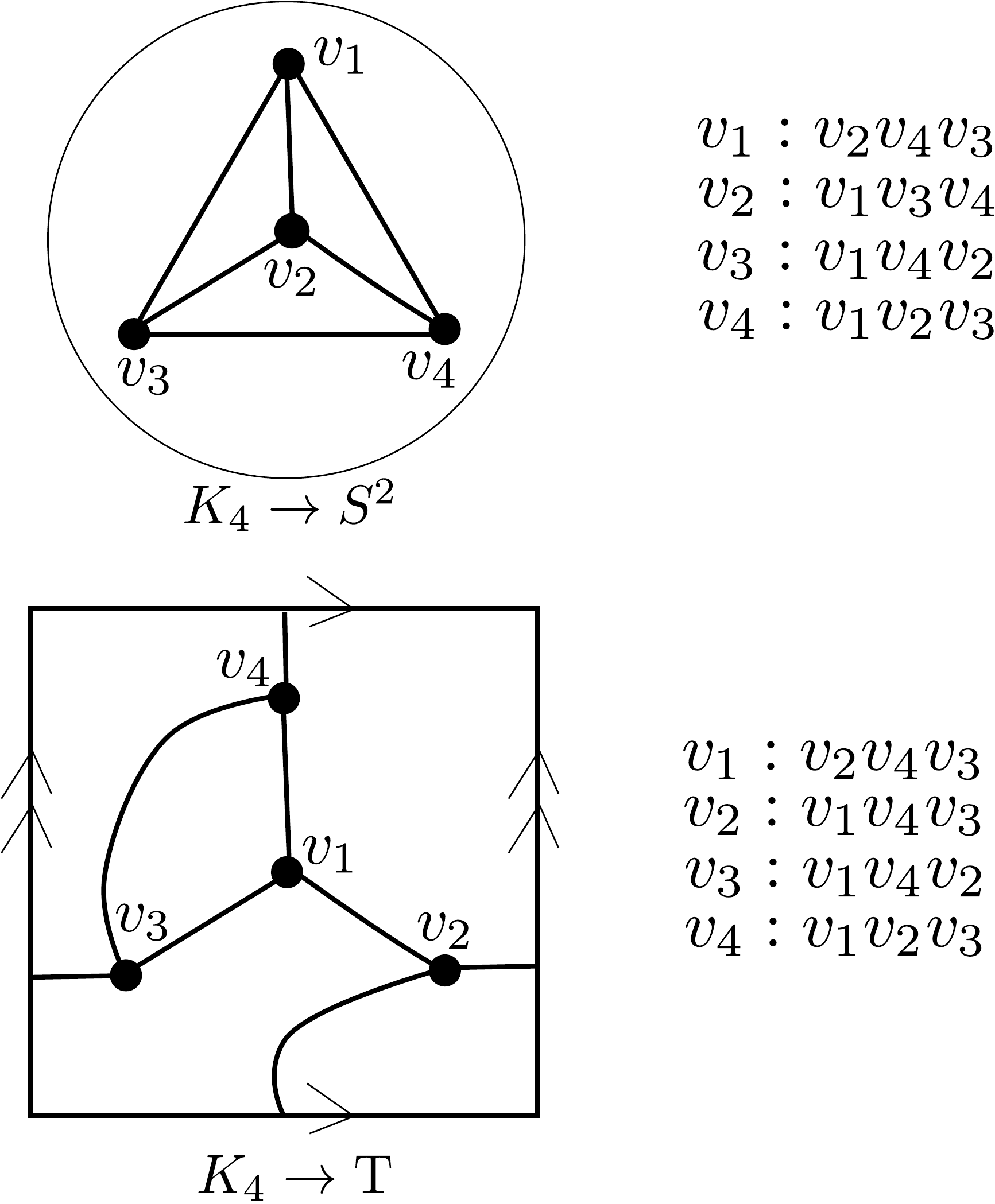}
\caption{An example of the impact of switching the order of two vertices in a rotation on a single vertex on a graph embedding: The switching ocurrs at the rotations on $v_2$.}\label{fig:RotationSwitched}
\end{center}
\end{figure}

We discuss here a matter left over from Section \ref{section:construction}. The reader will likely be curious about why it is that we do not add edges that will result in small angles being produced.  The answer is that we seek to maintain the accuracy and consistency of the appearances of vertices in all rotations in which they appear.  If vertices $a$ and $b$ appear in that order in a counterclockwise rotation about $v$, then they should appear in a reverse order in a counterclockwise rotation of $u$, which is a neighbor of $a$ and $b$ and far enough to the right of the projections of $a$ and $b$ in the tangent plane at $v$, as in Figure \ref{figure:RotationsReversed}.
However, in this case, the angle the angle between the edges joining the projections of $a$ and $b$ in the plane tangent at $v$ is small.  When one couples this fact with the notion that we are discarding points that are close to $v$ and to the tangent plane at $v$, the pitch of the tangent plane at $u$ is likely to be significantly different from that of the tangent plane at $v$.  Therefore the order in which $a$ and $b$ should appear in $u$ is reversed.  Note the impact of switching the order of $v_1$ and $v_4$ in the rotation on $v_2$ in Figure \ref{fig:RotationSwitched}. 

More rigorously, let $e_1$ and $e_2$ be unit vectors in $\mathbb{R}^n$, we consider the projection onto the plane spanned by two other unit vectors $t_1,t_2 \in \mathbb{R}^n$.  The coordinates of the projection of $e_1$ are $(x_1,y_1) = (e_1\cdot t_1,e_1\cdot t_2)$ and for $e_2$ the coordinates are $(x_2,y_2) = (e_2\cdot t_1, e_2 \cdot t_2)$, and we would like to insure that the orientation of these two 2-dimensional vectors is robust to small perturbations of $t_1$ and $t_2$.  The orientation is given by the sign of the determinant, namely the sign of 
\[ x_1 y_2-x_2 y_1 = (e_1 \cdot t_1)(e_2 \cdot t_2) - (e_2 \cdot t_1)(e_1 \cdot t_2) = t_1^\top (e_1 e_2^\top) t_2 - t_2^\top (e_1 e_2^\top) t_1 \]
where $e_1 e_2^\top$ is a rank-1 matrix whose only nonzero singular value is given by the trace $\textup{tr}(e_1e_2^\top) = \textup{tr}(e_2^\top e_1) = e_2\cdot e_1$.  Since all the other singular values are zero, it will be easiest to change the sign of the determinant when this non-zero singular value is close to zero.  Thus, by insuring that there are no small angles we ensure that small perturbations in the choice of tangent vectors $t_1,t_2$ will not result in different orientations of the projected vectors.


\subsection{Applying our methods to higher dimensions}\label{subsection:higherDimensions}
We should note that the method developed above for building a graph embedding as a way to understand a manifold could be applied manifolds of arbitrary dimension.  However, it only yields a full cell-complex for surfaces, and the rotation system can only be applied for surfaces.  In particular, our method could be used as a means of chosing a subset of a large data set that represents the topological information of the full data set.  Finding such a subset is an important problem in manifold learning and topological data analysis.  

On the other hand, for applications to surfaces, while the embedding dimension of the surface does not matter, the extrinsic curvature of the embedding determines the error in the tangent plane estimation.  So embeddings with high extrinsic curvature will require more data for this portion of the algorithm to work.

\subsection{Topics for further investigation}\label{subsection:topics}

The authors would like to put the following topics forward for investigation.
\begin{itemize}

\item{The authors used orthogonal projections into tangent planes as a way to produce a combinatorial representation of an embedded surface.  Is there some way to capture a better approximation by using some other kind of projection into a local approximation that is more form fitting than a tangent plane?}

\item{In light of Section \ref{section:smallAngles}, can we identify the offending rotation systems more easily then searching through all pairs of vertices?}

\item{Given a point cloud, how may we make better choices for the points at which we begin constructing the rotation system and the $\epsilon$ quantity that we chose?} 

\end{itemize}


\begin{thebibliography}{99}

\bibitem{Archdeacon}
D. Archdeacon. The medial graph and voltage-current duality. {\em Discrete Mathematics}, 104:111-141, 1992.

\bibitem{Berry}
T.~Berry and J.~Harlim. 
\newblock Iterated Diffusion Maps for Feature Identification.  
\newblock {\em Applied and Computational Harmonic Analysis}, 45(1):84--119, 2018.

\bibitem{Github}
\texttt{https://github.com/steveschluchter/rotation-systems-and-2-manifolds}

\bibitem{belkin2003laplacian}
M.~Belkin and P.~Niyogi.
\newblock Laplacian eigenmaps for dimensionality reduction and data
  representation.
\newblock {\em Neural computation}, 15(6):1373--1396, 2003.

\bibitem{diffusion}
R.~Coifman and S.~Lafon.
\newblock Diffusion maps.
\newblock {\em Appl. Comp. Harmonic Anal.}, 21:5--30, 2006.

\bibitem{SingerEstimate}
A.~Singer.
\newblock From graph to manifold {L}aplacian: The convergence rate.
\newblock {\em Applied and Computational Harmonic Analysis}, 21:128--134, 2006.

\bibitem{ghrist2008}
R.~Ghrist.
\newblock Barcodes: the persistent topology of data.
\newblock {\em Bulletin of the American Mathematical Society}, 45:61--75, 2008.

\bibitem{edels2010}
H.~Edelsbrunner and J.~Harer.
\newblock {\em Computational toplogy: an introduction}.
\newblock American Mathematical Soc., 2010.

\bibitem{carlsson2009topology}
G.~Carlsson.
\newblock Topology and data.
\newblock {\em Bulletin of the American Mathematical Society}, 46(2):255--308,
  2009.

\bibitem{GT}
J. Gross and T. Tucker, \textit{Topological graph theory}, Wiley, New York, 1987.

\bibitem{M} J. Munkres, \textit{Topology, Second Edition}, Prentice Hall, Upper Saddle River, 2000.





\end{thebibliography}
\end{document}